\documentclass[final]{elsarticle}

\usepackage{amssymb,latexsym}

\usepackage{amsmath}

\biboptions{sort&compress}

\journal{: \; Expositiones Mathematicae}

\begin{document}

\newtheorem{lema}{Lemma}
\newproof{prova}{Proof}
\newtheorem{teo}{Theorem}
\newtheorem{cor}{Corollary}

\begin{frontmatter}

\title{Some transcendence results from a harmless irrationality theorem}

\author{F. M. S. Lima}

\address{Institute of Physics, University of Bras\'{i}lia, P.O.~Box 04455, 70919-970, Bras\'{i}lia-DF, Brazil}


\ead{fabio@fis.unb.br}

\begin{abstract}
The arithmetic nature of values of some functions of a single variable, particularly, $\sin{z}$, $\cos{z}$, $\sinh{z}$, $\cosh{z}$, $e^z$, and $\ln{z}$, is a relevant topic in number theory. For instance, all those functions return transcendental values for all non-zero algebraic values of $z$ ($z \ne 1$ in the case of $\ln{z}$). On the other hand, not even an irrationality proof is known for some numbers like $\,e^e$, $\,\pi^e$, $\,\pi^\pi$, $\,\ln{\pi}$, $\,\pi + e\,$ and $\,\pi \, e$, though it is well-known that at least one of the last two numbers is irrational. In this note, I first derive a more general form of this last result, showing that at least one of the sum and product of any two transcendental numbers is transcendental. I then use this to show that, given any complex number $\,t \ne 0, 1/e$, at least two of the numbers $\,\ln{t}$, $\,t + e\,$ and $\,t \, e\,$ are transcendental. I also show that $\,\cosh{z}$, $\sinh{z}\,$ and $\,\tanh{z}\,$ return transcendental values for all $\,z = r \, \ln{t}$, $\,r \in \mathbb{Q}$, $r \ne 0$. Finally, I use a recent algebraic independence result by Nesterenko to show that, for all integer $\,n > 0$, $\,\ln{\pi}\,$ and $\,\sqrt{n} \,~ \pi\,$ are linearly independent over $\mathbb{Q}$.
\end{abstract}


\end{frontmatter}

\vspace{0.5cm}


As usual, let $\,\mathbb{Q}\,$ denote the set of all rational numbers, i.e.~the numbers which can be written as $\,p/\,q\,$, $p$ and $q$ being integers, $q \ne 0$. Also, let $\, \mathcal{A} \,$ denote the set of all algebraic numbers (over $\mathbb{Q}$), i.e.~the complex numbers $\,z\,$ which are the root of some (non-trivial) polynomial equation in $\,\mathbb{Q}[z]$. All other complex numbers, i.e.~$\,z \not \in \mathcal{A}\,$, are called \emph{transcendental} numbers. Of course, all rational numbers are algebraic, as they are roots of $\: q\,z -p = 0$, so all transcendental numbers are irrational. Though the existence of irrational numbers such as $\,\sqrt{2}\,$ remounts to the ancient Greeks, no example of a transcendental number was known at the beginning of the 19th century, which reflects the difficulty of proofing that a given number is transcendental.\footnote{The existence of `non-algebraic' numbers was conjectured by Euler in 1744, in his \textsl{Introduction to the analysis of the infinite}. There, he comments, without a proof, that ``the logarithms of (rational) numbers which are not powers of the base are neither rational nor (algebraic) irrational, so they should be called \emph{transcendental}.''}  A such proof appeared only in 1844, when Liouville showed that any number that has a rapidly converging sequence of distinct rational approximations must be transcendental~\cite{Liouv}.  In particular, he used his approximation theorem to show that $\,\sum_{k=0}^\infty{1/(2^{k!})}\,$ converges to a (real) transcendental number. However, it remained an important unsolved problem in pure mathematics to prove the transcendence of naturally occurring numbers.\footnote{From the pioneering work of Cantor on set theory in 1874, one knows that the set $\,\mathcal{A}\,$ is countable whereas the set $\,\mathbb{R}\,$ of all real numbers is uncountable, so \emph{`almost all' real numbers are transcendental}.} As a result, in 1873 Hermite proved that $\,e^r\,$ is transcendental for all rational $\,r \ne 0\,$ (in particular, $e\,$ is transcendental)~\cite{Hermite}. In 1882, Lindemann proved the following extension of Hermite's result~\cite{Lindemann}.

\begin{lema}[Hermite-Lindemann]
\label{lem:HL} \; The number $\,e^{\alpha}\,$ is transcendental for all algebraic $\,\alpha \ne 0$.
\end{lema}

This implies the transcendence of $\,\pi$, as follows from Euler's identity $\,e^{\,i\,\pi} = -1$, which is equivalent to the impossibility of squaring the circle with only ruler and compass, a problem that remained open by more than two thousand years. Note that Lemma~\ref{lem:HL} is equivalent to the transcendence of $\,\ln{\alpha}\,$ for all $\,\alpha \in \mathcal{A}\,$, $\,\alpha \ne 0, 1$. Hereafter, we are interpreting $\,\ln{z}\,$ as the principal value of this function, with the argument lying in the interval $(-\pi,\pi]$.

Based upon these first transcendence results, in 1885 Weierstrass succeeded in proofing a much more general result.

\begin{lema}[Lindemann-Weierstrass]
\label{lem:LW} \; Given a positive integer $\,n$, whenever $\,\alpha_0,\ldots,\alpha_n\,$ are distinct algebraic numbers, the numbers $\,e^{\alpha_0},\ldots,e^{\alpha_n}$ are linearly independent over $\mathcal{A}$. That is, for any $\,\beta_0,\ldots,\beta_n \in \mathcal{A}\,$ not all zero,
\begin{equation*}
\sum_{k=0}^n{\,\beta_k \, e^{\alpha_k}} \ne 0 .
\end{equation*}
\end{lema}

For a proof, see, e.g., Theorem 1.4 of Ref.~\cite{lt:Baker} or Theorem~1.8 of Ref.~\cite{lt:newBaker}.
\newline

As an immediate consequence, when one takes $\,\alpha_0 = 0\,$ and $\,\beta_0 \ne 0$, one concludes that

\begin{cor}
\label{cor0:LW} \; Given a positive integer $\,n$, being $\,\alpha_1, \ldots, \alpha_n\,$ distinct non-zero algebraic numbers and $\,\beta_1, \ldots, \beta_n \in \mathcal{A}\,$ not all zero, the sum $\,\sum_{k=1}^n{\beta_k \, e^{\alpha_k}}\,$ is a transcendental number.
\end{cor}

From the celebrated Euler's formula $\:e^{\pm \, i \theta} = \cos{\theta} \pm \, i \, \sin{\theta}$, it follows that $\,\cos{\theta} = (e^{i \theta} + e^{-i \theta})/2\,$ and $\,\sin{\theta} = (e^{i \theta} - e^{-i \theta})/(2 i)$. Analogously, the basic hyperbolic functions are defined as $\,\cosh{\theta} := (e^\theta + e^{-\theta})/2\,$ and $\,\sinh{\theta} := (e^\theta - e^{-\theta})/2$. From Corollary~\ref{cor0:LW}, it follows that

\begin{cor}
\label{cor1:LW} \; For any algebraic $\,\alpha \ne 0$, all numbers $\,\cos{\alpha}$, $\sin{\alpha}$, $\cosh{\alpha}$, and $\,\sinh{\alpha}$ are transcendental.
\end{cor}

The relevance of proofing the transcendence of powers and logarithms of algebraic numbers was acknowledged by Hilbert in his famous lecture ``Mathematical Problems'' in 1900, at the 2nd International Congress of Mathematicians~\cite{Hilbert}, being the content of his seventh problem, in which he questioned the arithmetic nature of $\,e^{i \pi z}\,$ for algebraic values of $\,z$.  Hilbert himself remarked that he expected this problem to be harder than finding a proof for the Riemann hypothesis, which remains open today.  Of course, for $\,z = r\,$ a rational it was already known that both $\cos{(r \pi)}$ and $\sin{(r \pi)}$ are algebraic, so $\,e^{i \pi r}\,$ is algebraic (from Euler's formula),\footnote{Moreover, it was proved by Lehmer in 1933 that, for rational $\,r=k/n$, $n>2$, the numbers $\,2 \cos{(2 \pi r)}\,$ and $\,2 \sin{(2 \pi r)}\,$ are \emph{algebraic integers} (i.e., roots of monic polynomial equations in $\mathbb{Z}[x]\,$). For details, see Ref.~\cite{Lehmer}.}  but the case of \emph{irrational} algebraic values of $\,z\,$ remained unsolved until 1934, when Gelfond and Schneider, working independently, showed that~\cite{Gelfond,Schneider}

\begin{lema}[Gelfond-Schneider]
\label{lem:GS} \;\, If $\:\alpha \ne 0,1\,$ and $\,\beta \not \in \mathbb{Q}\,$ are algebraic numbers, then any value of $\,\alpha^{\,\beta}$ is transcendental.
\end{lema}

For a proof, see e.g.~Theorem~2.1 of Ref.~\cite{lt:newBaker} or Theorem~10.1 (and Sec.~4 of Chap.~10) of Ref.~\cite{lt:Niven}.
\newline

This lemma promptly implies the transcendence of $\,2^{\sqrt{2}}\,$ and $\,e^\pi = i^{-2 i}$, two numbers mentioned by Hilbert in his lecture.\footnote{\label{ft:Hilbert}From Lemma~\ref{lem:GS}, the number $\,e^{i \pi \beta} = \left(e^{i \pi}\right)^\beta = (-1)^{\,\beta}$, which has the form $\,\alpha^{\,\beta}$, is transcendental for all $\,\beta \in \mathcal{A} \, \backslash \mathbb{Q}\,$, which solves Hilbert's 7th problem. Indeed, given $\,a, b \in \mathcal{A} \bigcap \mathbb{R}$, if either $\,a=0\,$ and $\,b \not \in \mathbb{Q}\,$ or $\,a \ne 0$, then the number $\,e^{(a +b i)\,\pi} = e^{(b -i a)\,i \pi} = (-1)^{b -i a}\,$ is transcendental.}  Note that Lemma~\ref{lem:GS} has a logarithmic version, namely

\begin{lema}[Log version]
\label{lem:logGS} \;\, The number $\:\log_{\beta}{\alpha} = \ln{\alpha}/\ln{\beta}\:$ is transcendental whenever $\,\alpha\,$ and $\,\beta\,$ are non-zero algebraic numbers, $\beta \ne 1$, and $\,\log_{\beta}{\alpha} \not \in \mathbb{Q}$.
\end{lema}

This form appears, e.g., in Theorem~10.2 of Ref.~\cite{lt:Niven}. It has a consequence for tangent arcs, as noted by Margolius in Ref.~\cite{Margolius}.

\begin{cor}[Margolius]
\label{cor1:GS} \; If $\,x\,$ is rational and $\: x \ne 0, \pm 1$, then the number
$\: \dfrac{\,\arctan{(x)}}{\pi} \:$
is transcendental.
\end{cor}

\begin{prova}
This result can be proved easily by writing $\,x = \tan{\theta}$, $x \in \mathbb{Q}$, $x \ne 0, \pm 1$. Then
\begin{equation}
\frac{\arctan{(x)}}{\pi} = \frac{\theta}{\pi} = \frac{1/i \cdot \ln{(z/|z|)}}{1/i \cdot \ln{(-1)}} = \frac{\ln{\left(\pm1 /\sqrt{1+x^2} + x \, i \, / \sqrt{1+x^2} \, \right)}}{\ln{(-1)}} \, ,
\label{eq:arctan}
\end{equation}
which comes on taking $\,z = \pm 1 +x\,i\,$ in $\,\ln{(z/|z|)} = i \, \theta$, which in turn comes from the exponential representation $\,z = |z| \, e^{i\,\theta}$. Clearly, the last expression in Eq.~\eqref{eq:arctan} is a ratio of two logs with algebraic arguments,\footnote{\label{fn:field}It follows from the fact that $\,\mathcal{A}\,$ is a \emph{field} that, given any $\,\alpha , \beta \in \mathcal{A}$, then all the numbers $\, \alpha \pm \beta$, $\alpha \, \beta$, and $\, \alpha / \beta\,$ ($\beta \ne 0$) are also algebraic (see, e.g., Sec.~6.6 and Theorem~6.12 of Ref.~\cite{Cameron}). More generally, given $\,r \in \mathbb{Q}\,$ and $\,\alpha \in \mathcal{A}$, $\alpha \ne 0$, if $\,z\,$ is any complex algebraic (respectively, transcendental) number then all numbers $\,z \pm \alpha$, $\alpha z$, $z/\alpha$, and $z^r$ are also algebraic (respectively, transcendental), the only exception being $\,z^0=1\,$ for $\,z \not \in \mathcal{A}$.} so Lemma~\ref{lem:logGS} applies and $\,\theta / \pi \,$ has to be either rational or transcendental. However, it is irrational because, being $\,r \in \mathbb{Q}$, $\,x = \tan{\theta} = \tan{(r\,\pi)}\,$ is rational only when $\,x = 0, \pm 1$, as proved by Niven in Corollary 3.12 of Ref.~\cite{lt:Niven}.\footnote{In fact, the irrationality of $\,\theta / \pi\,$ is nicely proved by Margolius by exploring the properties of sequences of primitive Pythagorean triples formed on writing $\,x = a/b$, $\,a$ and $b$ being distinct non-zero integers. See Theorem~3 of Ref.~\cite{Margolius}.}
\begin{flushright} $\Box$ \end{flushright}
\end{prova}

In particular, it follows that Plouffe's constant $\,\arctan{(\frac12)}/\pi\,$ is transcendental~\cite{Margolius}.  Let us extend Margolius' result to all basic trigonometric arcs. Hereafter, the word `trig' will stand for any of $\,\left\{\cos, \sin, \tan, \sec, \csc, \cot \right\}$.

\begin{teo}[Extension of Margolius' result]
\label{teo:extnsion} \; If $\,x\,$ is a real algebraic, then the number
$\; \dfrac{\mathrm{arctrig}{(x)}}{\pi} \;$
is either rational or transcendental.
\end{teo}

\begin{prova}
The proof is similar to the previous one, being enough to take $\,x = \mathrm{trig}{\,(\theta)}$, $x \in \mathcal{A} \bigcap \mathbb{R}$, and write
\begin{equation}
\frac{\mathrm{arctrig}{(x)}}{\pi} = \frac{\theta}{\pi} = \frac{\,\ln{(z/|z|)}}{\ln{(-1)}} \, ,
\label{eq:arctrig2}
\end{equation}
$z \ne 0$. Now, note that the choice of $\,z\,$ will change accordingly to the function represented by `trig'. For $\arccos{x}$, choose $\,z = x \pm \sqrt{1-x^2}\:i$. For $\mathrm{arcsec}{\,x}$, choose $\,z = \pm 1 \pm \sqrt{x^2-1} ~ i$. For $\arcsin{x}$, choose $\,z = \pm \sqrt{1-x^2} + x\,i$. For $\mathrm{arccsc}{\,x}$, choose $\,z = \pm \sqrt{x^2 -1} \pm i$. For $\arctan{x}$, choose $\,z = \pm 1 +x\,i$, as in the previous proof. For $\mathrm{arccot}{\,x}$, choose $\,z = x \pm i$. As the reader can easily check, in all these cases the ratio $\,z / |z|\,$ is an algebraic function of $\,x$, so it is an algebraic number for all real algebraic values of $\,x$. The last expression in Eq.~\eqref{eq:arctrig2} is then a ratio of two logs with algebraic arguments, which, from Lemma~\ref{lem:logGS}, we know to be either rational or transcendental.
\begin{flushright} $\Box$ \end{flushright}
\end{prova}

Conversely, if $\,x \in \mathbb{R}\,$ then it will be transcendental whenever $\,\mathrm{arctrig}{(x)}/\pi \in \mathcal{A} \, \backslash \,\mathbb{Q}\,$. This implies, e.g., the transcendence of $\,\mathrm{trig}\left(\sqrt{2} \; \pi \right)$.
\newline

The following extension of Lemma~\ref{lem:GS} was conjectured by Gelfond and proved by Baker in 1966, becoming the definitive result in this area.\footnote{Baker also gave a quantitative lower bound for these linear forms in logs, which had profound consequences for diophantine equations. This work won him a Fields medal in 1970.}

\begin{lema}[Baker]
\label{lem:Baker} \; Given non-zero algebraic numbers $\,\alpha_1, \ldots, \alpha_n\,$ such that $\,\ln{\alpha_1}, \ldots, \ln{\alpha_n}\,$ are linearly independent over $\,\mathbb{Q}$, then the numbers $\,1, \ln{\alpha_1}, \ldots, \ln{\alpha_n}\,$ are linearly independent over $\,\mathcal{A}$.  That is, for any $\,\beta_0,\ldots,\beta_n \in \mathcal{A}\,$ not all zero, we have
\begin{equation*}
\beta_0 + \sum_{k=1}^n{\,\beta_k \, \ln{\alpha_k}} \ne 0 \, .
\end{equation*}
\end{lema}

For a proof, see, e.g., Theorem 2.1 of Ref.~\cite{lt:Baker}.
\newline

This lemma has several interesting consequences.

\begin{cor}
\label{cor1:Baker} \; Given non-zero algebraic numbers $\,\alpha_1, \ldots, \alpha_n$, for any $\,\beta_1, \ldots, \beta_n \in \mathcal{A}\,$ the number $\,\beta_1 \, \ln{\alpha_1} + \ldots + \beta_n \, \ln{\alpha_n}\,$ is either null or transcendental. It is transcendental when $\,\ln{\alpha_1}, \ldots, \ln{\alpha_n}\,$ are linearly independent over $\,\mathbb{Q}\,$ and $\,\beta_1, \ldots, \beta_n\,$ are not all zero.
\end{cor}

For a proof, see, e.g., Theorem 2.2 of Ref.~\cite{lt:Baker}.

%

\begin{cor}
\label{cor3:Baker} \; Let $\,\alpha_1, \ldots, \alpha_n, \beta_0, \beta_1, \ldots, \beta_n\,$ be non-zero algebraic numbers. Then the number $\:e^{\,\beta_0} \, {\alpha_1}^{\beta_1} \ldots \, {\alpha_n}^{\beta_n}\,$ is transcendental.
\end{cor}

For a proof, see Theorem 2.3 of Ref.~\cite{lt:Baker}.

\begin{cor}
\label{cor4:Baker} \; For any algebraic numbers $\,\alpha_1, \ldots, \alpha_n\,$ other than $0$ or $1$, let $\,1,\beta_1, \ldots, \beta_n\,$ be algebraic numbers linearly independent over $\,\mathbb{Q}$. Then the number $\:{\alpha_1}^{\beta_1} \ldots \, {\alpha_n}^{\beta_n}\,$ is transcendental.
\end{cor}

For a proof, see Theorem 2.4 of Ref.~\cite{lt:Baker}.

\begin{cor}
\label{cor5:Baker} \; Let $\,\alpha \ne 0\,$ be an algebraic number. Then the number $\,e^{\alpha +\pi \, \beta}\,$ is transcendental for all algebraic values of $\:\beta$, without exceptions.
\end{cor}

For a proof, see Corollary 2 of Ref.~\cite{Lima2010}. Note that $\,e^{\alpha +\pi\beta}\,$ is transcendental even if $\,\alpha =0$, as long as $\,i \, \beta \not \in \mathbb{Q}\,$ (see Footnote~\ref{ft:Hilbert}).\footnote{Note also that Corollary~\ref{cor5:Baker} implies the transcendence of $\,(\alpha + \ln{\beta})/\pi\,$ for any non-zero $\,\alpha, \beta \in \mathcal{A}$.}
\newline

All this said, it is embarrassing that the numbers $\,\ln{\pi}$, $\pi + e\,$ and $\,\pi \, e\,$ are still not known to be transcendental. In fact, not even an irrationality proof is known, though it is easy to show that at least one of $\pi + e\,$ and $\,\pi \, e\,$ must be irrational. This is proved, e.g., in a nice survey on irrational numbers by Ross in Ref.~\cite{Ross}, but I include a short proof below for completeness.

Let us call \,\emph{quadratic}\, any algebraic number which is the root of a second-order polynomial equation in a single variable with rational coefficients. From the fact that $\,\pi\,$ is not a quadratic number (since it is not even an algebraic number, as commented just below our Lemma~\ref{lem:HL}) it follows that

\begin{lema}[Harmless irrationality]
\label{lem:naive} \; At least one of the numbers $\,\pi + e\,$ and $\,\pi \, e\,$ is irrational.
\end{lema}

\begin{prova}
Consider the quadratic equation $\,(x -\pi) \cdot (x -e) = 0$, whose roots are $\,\pi\,$ and $\,e$. By expanding the product, one has $\,x^2 -(\pi +e)\,x + \pi e = 0$. Assume, towards a contradiction, that both coefficients $\,\pi + e\,$ and $\,\pi \, e\,$ are rational numbers.  Then, our quadratic equation would have rational coefficients and both roots would be quadratic numbers. However, $\,\pi\,$ is \emph{not} a quadratic number.
\begin{flushright} $\Box$ \end{flushright}
\end{prova}

As the above proof depends only on the fact that $\,\pi\,$ is not a quadratic number and does not make use of any property of $\,e$, it is clear that Lemma~\ref{lem:naive} can be generalized.

\begin{lema}[General irrationality]
\label{lem:geral} \; Given any irrational number $\,u\,$ which is not quadratic and any complex number $\,v$, at least one of the numbers $\,u + v\,$ and $\,u \, v\,$ is irrational.
\end{lema}

\begin{prova}
The proof is identical to the previous one, being enough to substitute $\,\pi\,$ by $\,u\,$ and $\,e\,$ by $\,v$.
\begin{flushright} $\Box$ \end{flushright}
\end{prova}

In particular, this lemma applies when $\,u = t\,$ is a transcendental number, so at least one of the numbers $\,t + v\,$ and $\,t \, v\,$ is irrational.  Of course, for any algebraic $\,v \ne 0\:$ \emph{both} $\:t + v\,$ and $\,t \, v\,$ are \emph{transcendental} numbers,\footnote{\label{ft1}These basic transcendence rules are easily proved by contradiction.}  so the interesting case is when $\,v\,$ is also a transcendental number.  This leads us to the following transcendence result.

\begin{teo}[Transcendence of sums and products]
\label{teo:naiveT} \; Given two transcendental numbers $\,t_1\,$ and $\,t_2$, at least one of the numbers $\,t_1 + t_2\,$ and $\,t_1 \, t_2\,$ is transcendental.
\end{teo}

\begin{prova}
Given $\,t_1,t_2 \not \in \mathcal{A}$, consider the quadratic equation $\,(x -t_1) \, (x -t_2) = 0$, whose roots are $\,t_1\,$ and $\,t_2$. As it is equivalent to $\,x^2 -(t_1 +t_2)\,x + t_1 \, t_2 = 0$, assume, towards a contradiction, that both $\,s = t_1 + t_2\,$ and $\,p = t_1 \, t_2\,$ are algebraic numbers.  The equation then reads $\,x^2 -s \, x + p = 0$, so, by completing the square, one finds
\begin{eqnarray}
&&  x^2 -s x +\frac{s^2}{4} = \frac{s^2}{4} -p \nonumber \\
& \Longrightarrow & \left( x -\frac{s}{2} \right)^2 = \frac{s^2}{4} - p \, .
\end{eqnarray}
This implies that $\,(x -s/2)^2\,$ is algebraic (see Footnote~\ref{fn:field}), which is impossible because, being $\,x\,$ one of the roots $\,t_1\,$ and $\,t_2$, the number $\,x -s/2\,$ has to be transcendental (see Footnote~\ref{ft1}).
\begin{flushright} $\Box$ \end{flushright}
\end{prova}

This theorem implies, in particular, that at least one of $\,\pi+e\,$ and $\,\pi\,e\,$ is \emph{transcendental}. However, we are in a position to prove a stronger result.

\begin{teo}[Transcendence of two numbers]
\label{teo:1of3} \; Given any non-zero complex number $\,t \ne 1/e$, at least two of the numbers $\,t + e$, $\,t \, e\,$ and $\,\ln{t}\,$ are transcendental.
\end{teo}

\begin{prova}
If $\,t \ne 0\,$ is an algebraic number, then both $t+e$ and $t\,e$ are of course transcendental numbers, so let us restrict our attention to the transcendental values of $\,t$. If $\,\ln{t}\,$ is transcendental then we are done because we know, from Theorem~\ref{teo:naiveT}, that at least one of $\,t+e\,$ and $\,t\,e\,$ is transcendental. All that remains is to check whether $\,\ln{t} \in \mathcal{A}\,$ implies that both $\,t+e\,$ and $\,t\,e\,$ are transcendental numbers. Since $\,\ln{t} = \alpha \,\Longrightarrow \,t = e^\alpha$, then $\,t+e = e^\alpha+e\,$ is a transcendental number for all $\,\alpha \in \mathcal{A}$, $\alpha \ne 1$, according to Corollary~\ref{cor0:LW}. Note that $\,\alpha =1\,$ implies $\,t=e$, a case in which our theorem applies since $\,t+e = 2 \, e\,$ and $\,t \,e = e^2$ are both transcendental numbers. Also, since $\,\ln{t} \in \mathcal{A}\,$, then $\,1+\ln{t} = \ln{(e \, t)}\,$ is also algebraic, and then, according to Lemma~\ref{lem:HL}, the number $\,t \, e\,$ has to be transcendental for all $\,t\,$ such that $\,\ln{(t \, e)} \ne 0$, i.e. $t \ne 1/e\,$.
\begin{flushright} $\Box$ \end{flushright}
\end{prova}

In particular, this theorem implies that at least two of $\,\pi+e$, $\pi\,e\,$ and $\,\ln{\pi}\,$ are transcendental numbers.
\newline

Indeed, we can make suitable choices of $\,t_1\,$ and $\,t_2\,$ in Theorem~\ref{teo:naiveT} in order to obtain further transcendence results.

\begin{cor}
\label{cor:invT} \; For any transcendental number $\,t\,$ and algebraic numbers $\,\alpha\,$ and $\,\beta\,$ not both zero, the numbers $\,\alpha \, t + \beta / \, t\,$ and $\,t \, (\alpha - t)\,$ are both transcendental.
\end{cor}

\begin{prova}
For any $\,t \not \in \mathcal{A}\,$ and $\,\alpha,\beta \in \mathcal{A}\,$, not both zero, take $\,t_1 = \alpha \, t\,$ and $\,t_2 = \beta / \, t\,$  in Theorem~\ref{teo:naiveT}. If exactly one of $\,\alpha, \beta\,$ is null, then the proof is immediate. Otherwise, since $\,t_1 \, t_2 = \alpha \, \beta \in \mathcal{A}\,$ (see Footnote~\ref{fn:field}), then $\,t_1 +t_2 = \alpha \, t + \beta / \, t\,$ has to be a transcendental number. Finally, for any $\,\alpha \in \mathcal{A}$, take $\,t_1 = t\,$ and $\,t_2 = \alpha - t\,$ in Theorem~\ref{teo:naiveT}. Since $\,t_1 +t_2 = \alpha \in \mathcal{A}\,$, then the number $\,t_1 \, t_2 = t \, (\alpha -t)\,$ has to be transcendental.
\begin{flushright} $\Box$ \end{flushright}
\end{prova}

Given any transcendental number $t$, $\,t^r\,$ is transcendental for all rational $r$, $r \ne 0$.\footnote{This is easily proved by contradiction, writing $\,r = p/q$, $p\,$ and $\,q\,$ being non-zero integers.} What about $\,t^\alpha$, $\alpha\,$ being an irrational algebraic? On taking $\,t = e^\pi$, we know that $\,t^{\,i} \in \mathcal{A}$ whereas $t^{\sqrt{2}}$ is transcendental (see Footnote~\ref{ft:Hilbert}). The next theorem sheds some light on this question.

\begin{teo}[Existence of irrational exponent that retains transcendence]
\label{teo:algexpoente} \; For all transcendental number $\,t$, there is an irrational algebraic $\,\alpha\,$ such that $\,t^\alpha\,$ is also transcendental.
\end{teo}

\begin{prova}
Given any transcendental number $\,t$, assume (towards a contradiction) that $\,t^\alpha = \beta\,$ is algebraic for all $\,\alpha \in \mathcal{A} \, \backslash \mathbb{Q}\,$. From Corollary~\ref{cor:invT}, $\,t^r \, (\beta-t^r)\,$ is transcendental for all $\,r \in \mathbb{Q}$, $r \ne 0$, which means that $\,t^{r+\alpha} - t^{2 r} = t^{r+\alpha} \, (1-t^{r-\alpha})\,$ is transcendental. Clearly, $\alpha_{1,2} := r \pm \alpha\,$ is an irrational algebraic, so $\:t^{\alpha_1} \, (1-t^{\alpha_2}) = \beta_1 \, (1-\beta_2)\,$ should also be transcendental. However, it is the product of two algebraic numbers.
\begin{flushright} $\Box$ \end{flushright}
\end{prova}

Note that the similar proposition ``for all transcendental number $\,t$, there is an irrational algebraic $\,\alpha\,$ such that $\,t^\alpha\,$ is algebraic'' is \emph{false}, as follows on taking $\,t = e$ and using Lemma~\ref{lem:HL}.

Another consequence of Corollary~\ref{cor:invT} is

\begin{teo}[Linear independence of $\,1, \, \cosh{(r \ln{t})}$ and $\sinh{(r \ln{t})}$]
\label{teo:cosh} \; For any transcendental number $\,t\,$ and any rational $\,r \ne 0$, the numbers $\,1$, $\cosh{(r \ln{t})}$, and $\,\sinh{(r \ln{t})}\,$ are linearly independent over $\,\mathcal{A}\,$. In particular, $\cosh{(r \ln{t})}\,$ and $\,\sinh{(r \ln{t})}\,$ are transcendental numbers.
\end{teo}

\begin{prova}
Since $\,t^{\,r}\,$ is transcendental for any $\,t \not \in \mathcal{A}\,$ and any $\,r \in \mathbb{Q}$, $r \ne 0$, then, from Corollary~\ref{cor:invT}, one finds that, for any $\,\alpha,\beta \in \mathcal{A}$ not both zero:
\begin{eqnarray}
&& \alpha \, t^r + \frac{\beta}{t^r} = \alpha \, t^r + \beta \, t^{-r} \not \in \mathcal{A} \nonumber \\
& \Longrightarrow & \alpha \: e^{r \, \ln{t}} + \beta \: e^{- r \, \ln{t}} \not \in \mathcal{A} \nonumber \\
& \Longrightarrow & (\alpha + \beta) \: \cosh{(r \, \ln{t})} + (\alpha - \beta) \: \sinh{(r \, \ln{t})} \not \in \mathcal{A} \, .
\end{eqnarray}
Since $\alpha$ and $\beta$ are arbitrary algebraic numbers, then
\begin{equation}
\widetilde{\alpha} \: \cosh{(r \, \ln{t})} + \widetilde{\beta} \: \sinh{(r \, \ln{t})} \not \in \mathcal{A} \, ,
\end{equation}
where $\,\widetilde{\alpha} = \alpha +\beta\,$ and $\,\widetilde{\beta} = \alpha -\beta\,$ are also algebraic numbers (not both zero), so
\begin{equation}
\widetilde{\alpha} \: \cosh{(r \, \ln{t})} + \widetilde{\beta} \: \sinh{(r \, \ln{t})} \ne \gamma \, , \quad \forall \: \gamma \in \mathcal{A} \: .
\label{eq:boa}
\end{equation}
Therefore, $\, - \gamma + \widetilde{\alpha} \: \cosh{(r \, \ln{t})} + \widetilde{\beta} \: \sinh{(r \, \ln{t})} \ne 0$,
which shows that $\,1$, $\cosh{(r \ln{t})}\,$ and $\,\sinh{(r \ln{t})}\,$ are linearly independent over $\,\mathcal{A}\,$.

The transcendence of $\,\cosh{(r \ln{t})}\,$ follows on taking $\,\widetilde{\alpha} \ne 0 \,$ and $\,\widetilde{\beta} = 0\,$ in Eq.~\eqref{eq:boa}, whereas that of $\,\sinh{(r \ln{t})}\,$ follows on taking $\,\widetilde{\alpha} = 0\,$ and $\,\widetilde{\beta} \ne 0$.
\begin{flushright} $\Box$ \end{flushright}
\end{prova}

In addition, it is easy to prove the transcendence of $\,\tanh{(r \, \ln{t})}$.

\begin{teo}[Transcendence of $\,\tanh{(r \ln{t})}$]
\label{teo:tanh} \; For any transcendental number $\,t\,$ and any $\,r \in \mathbb{Q}$, $r \ne 0$, the number $\:\tanh{(r \ln{t})}\,$ is transcendental.
\end{teo}

\begin{prova}
For any transcendental number $\,t\,$ and any $\,r \in \mathbb{Q}$, $r \ne 0$, we have $\,\tanh{(r \ln{t})} := \sinh{(r \ln{t})}/\cosh{(r \ln{t})} = (t^r - t^{-r})/(t^r + t^{-r}) = (t^{2 r} -1)/(t^{2 r} +1) \ne 1$. Now, assume, towards a contradiction, that $\,\tanh{(r \ln{t})} = \alpha$, for some $\,\alpha \in \mathcal{A}$, $\alpha \ne 0, 1$. Then
\begin{eqnarray}
& & \frac{t^{2 r} -1}{t^{2 r} +1} = \alpha \nonumber \\
& \Longrightarrow & \,  t^{2 r} -1 = \alpha \, \left( t^{2 r} +1 \right) = \alpha \, t^{2 r} +\alpha \nonumber \\
& \Longrightarrow & \,  (1 -\alpha) \, t^{2 r} = \alpha +1 \nonumber \\
& \Longrightarrow & \,  t^{2 r} = \frac{1+\alpha}{1-\alpha} \, ,
\end{eqnarray}
which is impossible since the quotient of two algebraic numbers is also algebraic, whereas $\,t^{\,r} \not \in \mathcal{A}$.
\begin{flushright} $\Box$ \end{flushright}
\end{prova}

It follows, in particular, that $\,\cosh{(\ln{\pi})}$, $\sinh{(\ln{\pi})}\,$ and $\,\tanh{(\ln{\pi})}\,$ are transcendental numbers.  Interestingly, similar results can be derived for the basic trigonometric functions.

\begin{teo}[Linear independence of $\,1, \,\cos{(\beta \ln{\alpha})}\,$ and $\,\sin{(\beta \ln{\alpha})}\,$]
\label{teo:cos} \; For any algebraic numbers $\,\alpha, \beta$, $\alpha \ne 0, 1\,$ and $\: i \, \beta \not \in \mathbb{Q}$, the numbers $\,1$, $\cos{(\beta \ln{\alpha})}\,$ and $\,\sin{(\beta \ln{\alpha})}\,$ are linearly independent over $\,\mathcal{A}\,$. In particular, $\cos{(\beta \ln{\alpha})}\,$ and $\,\sin{(\beta \ln{\alpha})}\,$ are transcendental numbers.
\end{teo}

\begin{prova}
Since $\,i\,\beta\, \in \mathcal{A} \,\backslash \mathbb{Q}$, then, from Lemma~\ref{lem:GS}, $t = \alpha^{i \, \beta}\,$ is a transcendental number. From Corollary~\ref{cor:invT}, the sum $\,a \, t +b/t\,$ is also transcendental for all algebraic $a$ and $b$ not both zero, so
\begin{eqnarray}
&& a \, \alpha^{i \beta} + b \, \alpha^{- i \beta} = a \, e^{i \beta \, \ln{\alpha}} + b \, e^{- i \beta \, \ln{\alpha}} \not \in \mathcal{A} \nonumber \\
& \Longrightarrow & (a+b) \, \cos(\beta \, \ln{\alpha}) + i \, (a-b) \, \sin(\beta \, \ln{\alpha}) \not \in \mathcal{A} \, ,
\end{eqnarray}
which is equivalent to say that $\,\widetilde{a} \, \cos(\beta \, \ln{\alpha}) + \widetilde{b} \, \sin(\beta \, \ln{\alpha}) \ne c$, for all $\,c \in \mathcal{A}\,$. Therefore, for all $\,\widetilde{a},\widetilde{b} \in \mathcal{A}\,$ not both zero, $\widetilde{a} \, \cos(\beta \, \ln{\alpha}) + \widetilde{b} \, \sin(\beta \, \ln{\alpha}) - c \ne 0$, for all $\,c \in \mathcal{A}\,$. The transcendence of $\,\cos{(\beta \ln{\alpha})}\,$ follows by taking $\,\widetilde{a} \ne 0 , \, \widetilde{b} = 0\,$ and the transcendence of $\,\sin{(\beta \ln{\alpha})}\,$ follows by taking $\,\widetilde{a} = 0 , \, \widetilde{b} \ne 0$.
\begin{flushright} $\Box$ \end{flushright}
\end{prova}

\begin{teo}[Transcendence of $\,\tan{(\beta \ln{\alpha})}$]
\label{teo:tan} \; For any algebraic numbers $\,\alpha, \beta$, $\alpha \ne 0, 1\,$ and $\: i \beta \not \in \mathbb{Q}$, the number $\,\tan{(\beta \ln{\alpha})}\,$ is transcendental.
\end{teo}

\begin{prova}
Given non-zero algebraic numbers $\,\alpha, \beta$, $\alpha \ne 1$, assume, towards a contradiction, that $\,\tan{(\beta \ln{\alpha})} = \gamma \,$ for some $\,\gamma \in \mathcal{A}$. Then
\begin{eqnarray}
&& \gamma = \frac{\sin{(\beta \ln{\alpha})}}{\cos{(\beta \ln{\alpha})}} = \frac{1}{i} \, \frac{e^{i \beta \ln{\alpha}} -e^{-i \beta \ln{\alpha}}}{e^{i \beta \ln{\alpha}} +e^{i \beta \ln{\alpha}}} = \frac{1}{i} \, \frac{\alpha^{i \beta} -\alpha^{-i \beta}}{\alpha^{i \beta} +\alpha^{-i \beta}} \nonumber \\
& \Longrightarrow & \, i \, \gamma =  \frac{\, \alpha^{2 i \beta} -1 \,}{\alpha^{2 i \beta} +1} \, .
\end{eqnarray}
The last equality implies that $\,i \, \gamma \ne 1$. From Lemma~\ref{lem:GS}, $\,t = \alpha^{i \, \beta}\,$ is transcendental for all algebraic values of $\,\beta\,$ such that $\: i \, \beta \not \in \mathbb{Q}$. Then
\begin{eqnarray}
&& i \, \gamma = \widetilde{\gamma} = \frac{t^2 -1}{t^2 +1} \nonumber \\
& \Longrightarrow & \widetilde{\gamma} \; t^2 + \widetilde{\gamma} = t^2 -1 \nonumber \\
& \Longrightarrow & \left(1 -\widetilde{\gamma} \, \right) t^2 = 1 + \widetilde{\gamma} \nonumber \\
& \Longrightarrow & t^2 = \frac{\,1 +\widetilde{\gamma}\,}{1 -\widetilde{\gamma}} = \frac{\, 1 + \, i \, \gamma \,}{1 - \, i \, \gamma } \: ,
\end{eqnarray}
which should be algebraic since it is a quotient of two algebraic numbers. But this is impossible because, being $\,t\,$ a transcendental number, then $\,t^2\,$ is also transcendental.
\begin{flushright} $\Box$ \end{flushright}
\end{prova}

The last two theorems imply, for instance, that all $\,\mathrm{trig}(\ln{2})\,$ and $\,\mathrm{trig}\!\left( \sqrt{2} \; \pi \right)\,$ are transcendental numbers.
\newline

As a consequence of a recent work on modular functions by Nesterenko (1996)~\cite{Nesterenko}, we know that $\,\pi$ and $e^{\sqrt{n} \, \pi}$ are algebraically independent (over $\mathbb{Q}$) for all integers $\,n>0$. This has a consequence for the real number $\,\ln{\pi}$.

\begin{teo}[$\,\pi\,$ and $\,\ln{\pi}$]
\label{teo:lnpi} \; For any positive integer $\,n\,$ and any rational $\,r \ne 0$, the inequality
\begin{equation*}
q \, \ln{\pi} \ne \sqrt{n} ~ p \; \pi \, - \ln{r}
\end{equation*}
holds for all non-negative integers $\,p\,$ and $\,q$ (not both zero).
\end{teo}

\begin{prova}
Given $\,n \in \mathbb{Z}$, $n > 0$, since $\,\pi\,$ and $\,e^{\sqrt{n} \, \pi}\,$ are algebraically independent, then, for any non-zero integers $a$ and $b$ and non-negative integers $p$ and $q$, one has
\begin{eqnarray}
& & a \: \pi^{\,q} + b \left(e^{\sqrt{n} \: \pi}\right)^p \ne 0  \nonumber \\
\Longrightarrow & &  b \; e^{\sqrt{n} \: \pi \, p} \ne -a \: \pi^{\,q}  \nonumber \\
\Longrightarrow & &  e^{\sqrt{n} \: \pi \, p} \ne r \: \pi^{\,q} \, ,
\end{eqnarray}
where $\,r := -a/b\,$ is a non-zero rational. Since $\,e^{\sqrt{n} \: \pi \, p}\,$ and $\,\pi^q\,$ are both positive real numbers, then the last inequality obviously holds for any $\, r < 0$. For $\, r > 0$, we can take the logarithm on both sides, which yields
\begin{eqnarray}
& & \sqrt{n} \: p \: \pi \ne \ln{\left( r \, \pi^q \right)}  \nonumber \\
\Longrightarrow & &  \sqrt{n} \: p \: \pi \ne \ln{r} + q \, \ln{\pi}  \nonumber \\
\Longrightarrow & &  q \: \ln{\pi} \ne \sqrt{n} ~ p \; \pi - \ln{r} \, .
\end{eqnarray}
\begin{flushright} $\Box$ \end{flushright}
\end{prova}

On taking $\,r = 1$, one readily concludes that

\begin{cor}
\label{cor:pilnpi} \; For any positive integer $\,n$, the numbers $\,\ln{\pi}\,$ and $\,\sqrt{n} \; \pi\,$ are linearly independent over $\mathbb{Q}$.
\end{cor}

In other words, the number $\,\ln{\pi}\,$ is not a rational multiple of $\, \sqrt{n}\,~\pi$.\,\footnote{In particular, $\ln{\pi}\,$ is not a rational multiple of $\pi$.}
\newline

I hope the results put forward here in this paper can be useful for those who are studying the arithmetic nature of numbers such as $\,\ln{\pi}$, $\pi \, e\,$, $\pi + e$, $e^e$, $\pi^\pi\,$ and other related numbers.


\section*{Acknowledgments}
The author thanks Mr.~Bruno S.~S.~Lima and Mrs.~Claudionara de Carvalho for all helpful discussions on classical transcendence results.


\end{document}